\documentstyle[12pt]{article}
\begin{document}
\newtheorem{Def}{Definition}
\newtheorem{thm}{Theorem}
\newtheorem{lem}{Lemma}
\newtheorem{rem}{Remark}
\newtheorem{prop}{Proposition}
\newtheorem{cor}{Corollary}
\title
{A general Liouville type theorem
 for some conformally invariant fully nonlinear equations } 
\author{\ Aobing Li
\ \ \ \ \& \ \ \ YanYan Li\thanks{Partially supported by
NSF Grant DMS-0100819.}
\\ Department of Mathematics\\ Rutgers University\\
110 Frelinghuysen Rd.\\
Piscataway, NJ 08854
}
\date{}
\maketitle
\newcommand{\Bbb}{ }
%\setcounter{section}{0}
%\section{}

Various Liouville type theorems for conformally invariant equations
have been obtained by Obata (\cite{O}),  Gidas, Ni and Nirenberg
(\cite{GNN}),
Caffarelli, Gidas and Spruck (\cite{CGS}), Viaclovsky
{(\cite{V3} and \cite{V2}), Chang, Gursky and Yang
(\cite{CGY2} and \cite{CGY3}), and Li and Li (\cite{LL0}, \cite{LL}
and \cite{LL2}).  See e. g. theorem 1.3 and remark 1.6 in \cite{LL}
where  these results (except for the one in \cite{LL2}) are stated
more precisely.

In this paper we give a general Liouville type theorem for
conformally invariant fully nonlinear equations.  This
extends the above mentioned Liouville type theorems.

For $n\ge 3$, 
let ${\cal S}^{n\times n}$ be the set of  $n\times n$ real symmetric
matrices,   
${\cal S}^{n\times n}_+\subset {\cal S}^{n\times n}$ be the set
of positive definite matrices, and let $O(n)$ be 
the set of $n\times n$ real orthogonal matrices.

For a positive $C^2$ function $u$, let
$$
A^u:= -\frac{2}{n-2}u^{  -\frac {n+2}{n-2} }
\nabla^2u+ \frac{2n}{(n-2)^2}u^ { -\frac {2n}{n-2} }
\nabla u\otimes\nabla u-\frac{2}{(n-2)^2} u^ { -\frac {2n}{n-2} }
|\nabla u|^2I,
$$
where  $I$ is the $n\times n$ identity matrix.

Let $U\subset{\cal S}^{n\times n}$ be an open set satisfying 
\begin{equation}
O^{-1}UO=U,\quad\forall~O\in O(n)
\label{1}
\end{equation}
and
\begin{equation}
U\cap\{M+tN|~0<t<\infty\}~~\mbox{is convex}\quad\mbox{for}~
\forall~M\in{\cal S}^{n\times n}, N\in{\cal S}^{n\times n}_+.
\label{2}
\end{equation}
Let $F\in C^1 (U)$ satisfy 
\begin{equation}
F(O^{-1}MO)=F(M),\quad\forall M\in U, O\in O(n)
\label{3}
\end{equation}
and 
\begin{equation}
(F_{ij}(M))>0,\quad\forall M\in U,
\label{4}
\end{equation}
where $F_{ij}(M)=\frac{\partial F}{\partial M_{ij}}(M)$.\newline

\begin{thm}
For $n\ge 3$, let $U\subset{\cal S}^{n\times n}$ be open and satisfy 
(\ref{1}) and (\ref{2}), and let $F\in C^1(U)$ satisfy (\ref{3}) and (\ref{4}). 
Assume that $u\in C^2(\Bbb{R}^n)$ is a positive function satisfying 
\begin{equation}
F(A^u)=1,~~~A^u\in U,\quad\mbox{on}~~\Bbb{R}^n,
\label{5}
\end{equation}
and 
\begin{equation}
\Delta u\le 0,\quad\mbox{on}~~\Bbb{R}^n.
\label{6}
\end{equation}
Then for some $\bar x\in \Bbb{R}^n$ and some constants $a>0$ and 
$b\ge 0$ satisfying $2b^2a^{-2}I\in U$ and $F(2b^2a^{-2}I)=1$
\begin{equation}
u(x)\equiv\Big( \frac{a}{1+b^2|x-\bar x|^2} \Big)^{\frac{n-2}{2}},
\quad\forall x\in \Bbb{R}^n.
\label{7}
\end{equation}
\label{theorem1}
\end{thm}
\begin{rem}
If $U$ has the property that 
\begin{equation}
Trace(M):=\sum_{i=1}^n M_{ii}\ge 0,\quad\forall~M\in U,
\label{8}
\end{equation}
then any positive solution $u$ of (\ref{5}) automatically satisfies (\ref{6}).
\label{remark1}
\end{rem}
\begin{rem}
When $b=0$ in (\ref{7}), then $u\equiv$Constant, $A^u\equiv 0$, 
 $0\in U$, and $F(0)=1$.
\label{remark1a}
\end{rem}

Let $B_R(x)\subset\Bbb{R}^n$ denote the ball of radius $R$ centered at 
$x$, and let $B_R=B_R(0)$.
\begin{lem}
For $n\ge 1,~R>0$, let $\xi\in C^2(B_R\setminus\{0\})$ satisfy 
\begin{equation}
\Delta\xi\ge 0\quad\mbox{in}~B_R\setminus\{0\},
\label{9}
\end{equation}
and 
\begin{equation}
\inf\limits_{B_R\setminus\{0\}}\xi >-\infty.
\label{10}
\end{equation}
Assume that there exist $\eta,~\zeta\in C^1(B_R)$ satisfying 
\begin{equation}
\Delta \eta\ge 0,~~\Delta\zeta\ge 0,\quad\mbox{in}~B_R~~
\mbox{in the distribution sense},
\label{11}
\end{equation}
\begin{equation}
\eta(0)=\zeta (0),
\label{12}
\end{equation}
\begin{equation}
\nabla \eta(0)\neq\nabla\zeta (0),
\label{13}
\end{equation}
and
\begin{equation}
\xi\le\eta,\quad\xi\le\zeta,\quad\mbox{in}~B_R\setminus\{0\}.
\label{14}
\end{equation}
Then
\begin{equation}
\limsup\limits_{x\to 0}\xi (x)<\eta (0).
\label{15}
\end{equation}
\label{lemma1}
\end{lem}
\begin{rem}
If we further assume that $\eta,\zeta\in C^2(B_R)$, then hypothesis 
(\ref{10}) is not needed in Lemma~\ref{lemma1}. This can be deduced 
easily from Lemma~\ref{lemma2} by letting $\xi=-u$, $\eta=-w$ and 
$\zeta=-v$.
\label{remmark2}
\end{rem}
{\bf Proof of Lemma~\ref{lemma1}.}\ Replacing $\xi$, $\eta$ and $\zeta$ by 
$$
\tilde\xi(x)=\xi(x)-\nabla\eta(0)\cdot x
+|\nabla\eta(0)|R+1-\inf\limits_{B_R\setminus\{0\}}\xi,
$$ 
$$
\tilde \eta(x)=\eta(x)-\nabla\eta(0)\cdot x
+|\nabla\eta(0)|R+1
-\inf\limits_{B_R\setminus\{0\}}\xi,
$$ 
and 
$$
\tilde \zeta(x)=\zeta (x)-\nabla\eta(0)\cdot x
+|\nabla\eta(0)|R+1
-\inf\limits_{B_R\setminus\{0\}}\xi
$$
respectively, we may further assume that 
\begin{equation}
\nabla\eta (0)=0,
\label{16}
\end{equation}
\begin{equation}
\xi\ge 1\quad\mbox{in}~B_R\setminus\{0\}.
\label{17}
\end{equation}
Without loss of generality, we may assume that $R=1$. By (\ref{13}) and (\ref{16}), 
$\nabla\zeta(0)\neq 0$. After making a rotation, we may assume that
\begin{equation}
\frac{\partial\zeta}{\partial x_1}(0)=-|\nabla\zeta (0)|<0.
\label{18}
\end{equation}
Since $\xi\in L_{loc}^{\infty}(B_1)$ and $\Delta\xi\ge 0$ in $B_1\setminus\{0\}$, 
we know that $\Delta\xi\ge 0$ in $B_1$ in the distribution sense. Consequently,
\begin{equation}
\xi(y)\le\frac{1}{|B_r(y)|}\int_{B_r(y)}\xi,\quad\forall~0<|y|<1,~\forall~0<r<1-|y|.
\label{19}
\end{equation}
Since $\zeta\in C^1(B_1)$ satisfies (\ref{18}), there exists $0<\delta<\frac 12$ 
such that for $\forall~e\in\Bbb{R}^n$ with 
$|e|=1$ and $e\cdot e_1\ge 1-\delta$, we have 
\begin{equation}
\nabla\zeta (x)\cdot e<-\delta,\quad\forall~|x|<\delta.
\label{20}
\end{equation}
where $e_1=(1,0,\cdots,0)$.\newline
Let $S_{\delta}:=\{x\in\Bbb{R}^n\setminus\{0\}|~\frac{x}{|x|}\cdot e_1>1-\delta\}$. 
Now we fix the value of $\delta$. In the following, we will choose small positive 
numbers $r$ and $t$ satisfying $0<t<\frac{r}{10}<\frac{\delta}{40}$, and we will 
show that for some positive constant $c$,
depending only  on $\delta,~n$ and $r$, 
we have 
\begin{equation}
\xi(y)\le\eta(0)-c,\quad\forall~0<|y|<t.
\label{21}
\end{equation}
For $0<|y|<t$, we have, by using (\ref{14}), (\ref{16}), (\ref{19}), that
\begin{eqnarray*}
\xi(y)&\le&\frac{1}{|B_r(y)|}\int_{B_r(y)}\xi\le\frac{1}{|B_r(y)|}
\{\int_{B_r(y)\setminus S_{\delta}}\eta+\int_{B_r(y)\cap S_{\delta}}\zeta\}\\
&&=\frac{1}{|B_r(y)|}\{\int_{B_r(y)\setminus S_{\delta}}(\eta (0)+o(r))
+\int_{B_r(y)\cap S_{\delta}}\zeta\},
\end{eqnarray*}
where $o(r)$ satisfying $\lim\limits_{r\to 0}\frac{|o(r)|}{r}=0$.\newline
First recall that $|y|<t<\frac{r}{10}$, 
\[
\int_{B_r(y)\setminus S_{\delta}}(\eta (0)+o(r))
=\eta (0)|B_r(y)\setminus S_{\delta}|+o(r^{n+1}).
\]
Next recall that $\zeta$ satisfies (\ref{20}), $1\le\xi\le\zeta$ and 
$|y|<t<\frac{r}{10}$,
\begin{eqnarray*}
\int_{B_r(y)\cap S_{\delta}}\zeta&\le&\int_{B_{r+|y|\cap S_{\delta}}}\zeta 
=\int_0^{r+|y|}(\int_{\partial B_s\cap S_{\delta}}\zeta)~ds\\
&\le&\int_0^{r+|y|}\Big( \int_{\partial B_s\cap S_{\delta}}(\zeta(0)-\delta s) \Big)~ds\\
&=&\zeta(0)|B_{r+|y|}\cap S_{\delta}|-\frac{\delta }{n+1}|\partial B_1\cap S_{\delta}|(r+|y|)^{n+1}\\
&\le&\zeta(0)|B_{r+2|y|}(y)\cap S_{\delta}|-\frac{\delta }{n+1}|\partial B_1\cap S_{\delta}|r^{n+1}
\end{eqnarray*}
Since $\zeta(0)=\eta(0)$, we deduce from the above that 
\begin{eqnarray*}
\xi(y)&\le&\frac{1}{|B_r(y)|}\{\eta (0)|B_r(y)\setminus S_{\delta}|+o(r^{n+1})+
\eta(0)|B_{r+2|y|}(y)\cap S_{\delta}|\\
&&-\frac{\delta}{n+1}|\partial B_1\cap S_{\delta}|r^{n+1}\}\\
&=&\eta(0)+\frac{1}{|B_r(y)|}\{\eta(0)|(B_{r+2|y|}(y)\setminus B_r(y))\cap S_{\delta}|\\
&&+ o(r^{n+1})-\frac{\delta}{n+1}|\partial B_1\cap S_{\delta}|r^{n+1} \}
\end{eqnarray*}
Now fix some small $r$ satisfying $0<r<\frac{\delta}{4}$ and 
$o(r^{n+1})-\frac{\delta}{2(n+1)}|\partial B_1\cap S_{\delta}|r^{n+1}\le 0$. 
Since 
\[
|B_{r+2|y|}(y)\setminus B_r(y)|\le C(n) r^{n-1}|y|,
\] 
we can fix a smaller $t$ satisfying $0<t<\frac{r}{10}$ such that 
\[
C(n)\eta(0)r^{n-1}t-\frac{\delta}{4(n+1)}|\partial B_1\cap S_{\delta}|r^{n+1}\le 0.
\] 
With these choices of $r$ and $t$, we have 
\[
\xi (y)\le \eta (0)-\frac{\delta}{4(n+1)}|\partial B_1\cap S_{\delta}|r^{n+1}.
\] 
Estimate (\ref{15}) follows from the above. 
Lemma~\ref{lemma1} is established.

\vskip 5pt
\hfill $\Box$
\vskip 5pt

\begin{lem}
For $n\ge 1$, $R>0$, let $u\in C^2(B_R\setminus\{0\})$ satisfy $\Delta u\le 0$ in 
$B_R\setminus\{0\}$. Assume that there exist $w,~v\in C^2(B_R)$ satisfying 
\[
\Delta w\le 0,~~\Delta v\le 0\quad\mbox{in}~B_R,
\]
\[
w(0)=v(0),\quad \nabla w(0)\neq\nabla v(0),
\]
and 
\[
u\ge w,\quad u\ge v,\quad\mbox{in}~B_R\setminus\{0\}.
\]
Then 
$$\liminf\limits_{x\to 0}u(x)>w(0).
$$
\label{lemma2}
\end{lem}
{\bf Proof of Lemma~\ref{lemma2}.}\ By adding a large constnat to $u,~w$ and $v$, 
we may assume that 
\[
v\ge 1,\quad w\ge 1,\quad\mbox{in}~B_{\frac R2}\setminus\{0\}.
\]
Let $\xi=\frac 1u$, $\eta=\frac 1w$ and $\zeta=\frac 1v$. 
Since $\Delta u\le 0$ in 
$B_R\setminus\{0\}$, a straight forward calculation yields 
\[
\Delta \xi=-u^{-2}\Delta u+2u^{-3}|\nabla u|^2\ge 0,\quad\mbox{in}
~B_{\frac R2}\setminus\{0\}.
\]
Similarly, we have 
\[
\Delta \eta\ge 0,\quad\Delta\zeta\ge 0,\quad\mbox{in}~B_{\frac R2}.
\]
Clearly, $\eta(0)=\zeta(0)$, $\nabla\eta(0)\neq\nabla\zeta(0)$ and $\xi>0$ in 
$B_{\frac R2}\setminus\{0\}$. It follows from Lemma~\ref{lemma1} that 
$\limsup\limits_{x\to 0}\xi(x)<\eta(0)$, i.e., $\liminf\limits_{x\to 0}u(x)
>w(0)$. 
Lemma~\ref{lemma2} is established.

\vskip 5pt
\hfill $\Box$
\vskip 5pt

\begin{prop}
For $n\ge 3$, let $U\subset {\cal S}^{n\times n}$ be open and satisfy (\ref{1}) and (\ref{2}), 
and let $F\in C^1(U)$ satisfy (\ref{3}) and (\ref{4}). Assume that 
$u\in C^2(\Bbb{R}^n\setminus\{0\})$ is a positive function satisfying
\begin{equation}
F(A^u)=1,~~~A^u\in U,\quad\mbox{in}~\Bbb{R}^n\setminus\{0\},
\label{22}
\end{equation}
\begin{equation}
\Delta u\le 0,\quad\mbox{in}~\Bbb{R}^n\setminus\{0\},
\label{23}
\end{equation}
and 
\begin{equation}
u_{0,1}~~\mbox{can be extended to a}~C^2~\mbox{function near the origin},
\label{24}
\end{equation}
where $u_{0,1}(x):=\frac{1}{|x|^{n-2}}u(\frac{x}{|x|^2})$.\newline
We further assume that
there exist some constant $\delta>0$ and $v\in C^2(B_{\delta})$ 
such that
\begin{equation}
\Delta v\le 0\quad\mbox{in}~B_{\delta}.
\label{25}
\end{equation}
\begin{equation}
\nabla v(0)=0,
\label{26}
\end{equation}
\begin{equation}
\label{proposition1}
u-v\ge 0\quad\mbox{in}~B_{\delta}\setminus\{0\},
\label{27}
\end{equation}
\begin{equation}
\liminf\limits_{x\to 0}(u-v)(x)=0.
\label{27a}
\end{equation}
Then $u$ is radially symmetric, i.e,
\begin{equation}
u(x)=u(y),\quad\forall~|x|=|y|>0.
\label{27b}
\end{equation}
Moreover $u'(r)<0$ for $\forall~r>0$, where we have used $u(r)$ to 
denoted the radially symmetric function $u$.
\label{prop1}
\end{prop}

\begin{lem}
Let $u\in C^0(B_2\setminus\{0\})$ satisfy 
\[
\Delta u\le 0\quad\mbox{in}~B_2\setminus\{0\}~\mbox{in the distribution sense},
\]
and $\inf\limits_{B_2\setminus\{0\}}u>-\infty$. Then 
\[
u\ge \min\limits_{\partial B_1}u\quad\mbox{on}~B_1\setminus\{0\}.
\]
\label{lemma3}
\end{lem}
{\bf Proof of Lemma~\ref{lemma3}.}\ For $\epsilon>0$, consider 
$v_{\epsilon}(x):=\epsilon (1-\frac{1}{|x|^{n-2}})+\min\limits_{\partial B_1} u$. 
Then 
\[
\Delta (v_{\epsilon}-u)\ge 0~~~~\mbox{in}~B_1\setminus\{0\},\quad
(v_{\epsilon}-u)\le 0~~~~\mbox{on}~\partial B_1.
\]
Since $\limsup\limits_{x\to 0}(v_{\epsilon}(x)-u(x))=-\infty$, we 
deduce from the maximum principle that   
\[
v_{\epsilon}-u\le 0\quad\mbox{on}~B_1\setminus\{0\}.
\]
Fix any $x$ in $B_1\setminus\{0\}$, and send $\epsilon\to 0$, we have 
$u(x)\ge \min\limits_{\partial B_1}u$. Lemma~\ref{lemma3} is established.

\vskip 5pt
\hfill $\Box$
\vskip 5pt

\noindent
{\bf Proof of Proposition~\ref{prop1}.}\ By the positivity of $u$
and by (\ref{23}), we have 
$u_{0,1}>0$ and $\Delta u_{0,1}\le 0$ on $\Bbb{R}^n\setminus\{0\}$. By Lemma~\ref{lemma3}, 
\begin{equation}
\inf\limits_{B_1\setminus\{0\}}u>0,\quad\min\limits_{B_1}u_{0,1}>0.
\label{28}
\end{equation}
If $u$ can be extended to a $C^1$ function near the origin, then, by theorem~1.2 in 
\cite{LL}, 
$u$ is of the form (\ref{7}) for some $\bar x\in\Bbb{R}^n$ and some positive constants 
$a$ and $b$. By (\ref{27}), (\ref{27a}) and
(\ref{26}),  $\nabla u(0)=0$, and therefore $\bar x=0$. 
Proposition~\ref{prop1} is proved in this case. 
{\it In the rest of the proof of 
Proposition~\ref{prop1}, we always assume that $u$ can not be extended 
to a $C^1$ function near the origin.}

By (\ref{28}) and the repeatedly
used  arguments in  \cite{LZ}, \cite{LL} and \cite{LL2},
 we can prove that 
$\forall~x\in\Bbb{R}^n\setminus\{0\}$, there exists $\lambda_0(x)>0$ such that
\[
u_{x,\lambda}(y):=(\frac{\lambda}{|y-x|})^{n-2}u(x+\frac{\lambda^2(y-x)}{|y-x|^2})\le u(y),
\quad\forall~0<\lambda<\lambda_0(x),~|y-x|\ge \lambda,~y\neq 0.
\]
Set 
$$\bar\lambda(x)=\sup\{0<\mu<|x|~|~u_{x,\lambda}(y)\le u(y),~\forall~|y-x|\ge\lambda,
~y\neq 0,~0<\lambda\le\mu\}.
$$

We distinguish into two cases.\newline
{\bf Case 1.}\ $\exists~\bar x\in\Bbb{R}^n\setminus\{0\}$ such that $\bar\lambda(\bar x)<|\bar x|$.\newline
{\bf Case 2.}\ $\bar\lambda(x)=|x|$ for $\forall~x\in\Bbb{R}^n\setminus\{0\}$.\newline

{\bf  In Case 1}, we have 
\begin{equation}
u_{\bar x,\lambda}(y) \le u(y),
\quad\forall~0<\lambda<\bar \lambda(\bar x),~|y-\bar x|\ge \lambda,~y\neq 0.
\label{28a}
\end{equation}
After a rotation, we may assume that $\bar x=\bar x_1e_1$ with $\bar x_1>0$.
\begin{lem}
$\nabla u_{\bar x,\bar\lambda (\bar x)}(0)\neq 0$.
\label{lemma4}
\end{lem}
{\bf Proof of Lemma~\ref{lemma4}.}\ Suppose the contrary, 
\begin{equation}
\nabla u_{\bar x,\bar\lambda(\bar x)}(0)=0.
\label{29}
\end{equation}
A direct calculation yields that
\[
\partial_{y_1}u_{\bar x,\bar\lambda(\bar x)}(0)=(n-2)\bar\lambda(\bar x)^{n-2}|\bar x|^{1-n}
u((1-(\frac {\bar\lambda(\bar x)}{|\bar x|})^2)\bar x)-\bar\lambda(\bar x)^n|\bar x|^{-n}
\partial_1 u((1-(\frac {\bar\lambda(\bar x)}{|\bar x|})^2)\bar x).
\]
By (\ref{29}),
\begin{equation}
(n-2)|\bar x| u((1-(\frac {\bar\lambda(\bar x)}{|\bar x|})^2)\bar x)=
\bar\lambda(\bar x)^2\partial_1 u((1-(\frac {\bar\lambda(\bar x)}{|\bar x|})^2)\bar x).
\label{30}
\end{equation}
Consider $w(s):=u(\bar x-s\frac{\bar x}{|\bar x|})$ for $s>0$. By (\ref{28a}) with 
$y=\bar x-s\frac{\bar x}{|\bar x|}$, 
\[
(\frac{\lambda}{s})^{n-2}w(\frac{\lambda^2}{s})\le w(s),
\quad\forall~\lambda\le s<|\bar x|,~\forall~0<\lambda\le\bar\lambda(\bar x).
\]
It follows (with $t=\frac{\lambda^2}{s}$) that
 $t^{\frac{n-2}{2}}w(t)\le s^{\frac{n-2}{2}}w(s)$ 
 $\forall~0<t\le s\le\bar\lambda(\bar x)$, and therefore 
(note that $\frac{\bar\lambda(\bar x)^2}{|\bar x|}<\bar\lambda(\bar x)$) 
\[
\frac{d}{ds}(s^{\frac{n-2}{2}}w(s))|_{s=\frac{\bar\lambda(\bar x)^2}{|\bar x|}}\ge 0,
\]
i.e.,
\begin{equation}
\frac{n-2}{2}u((1-(\frac{\bar\lambda(\bar x)}{|\bar x|})^2)\bar x)\ge 
\frac{\bar\lambda(\bar x)^2}{|\bar x|}\partial_1 u((1-(\frac{\bar\lambda(\bar x)}{|\bar x|})^2)\bar x).
\label{31}
\end{equation}
By (\ref{30}) and (\ref{31}), $\frac{n-2}{2}u((1-(\frac{\bar\lambda(\bar x)}{|\bar x|})^2)\bar x)\ge 
(n-2)u((1-(\frac{\bar\lambda(\bar x)}{|\bar x|})^2)\bar x)$. This is a contradiction, 
since $u((1-(\frac{\bar\lambda(\bar x)}{|\bar x|})^2)\bar x)$ and $n-2>0$.
 Lemma~\ref{lemma4} is established.

\vskip 5pt
\hfill $\Box$
\vskip 5pt

Since $\bar\lambda(\bar x)<|\bar x|$, we have, by (\ref{28a}), that $u_{\bar x,\bar\lambda(\bar x)}\le u$ 
in an open neighborhood of the origin. Since $u$ is a $C^2$ superharmonic function in 
$\Bbb{R}^n\setminus\{0\}$, $u_{\bar x,\bar\lambda(\bar x)}(\bar x)$ is a superharmonic function 
in an open neighborhood of the origin. We first show that
\begin{equation}
\liminf\limits_{|y|\to \infty}|y|^{n-2}(u-u_{\bar x,\bar\lambda(\bar x)})(y)>0.
\label{32}
\end{equation}
Indeed, let $\xi (x)=\frac{1}{|x|^{n-2}}u(\frac{x}{|x|^2})$ and 
$\eta(x)=\frac{1}{|x|^{n-2}}u_{\bar x,\bar\lambda(\bar x)}(\frac{x}{|x|^2})$. 
By the hypothesis on $u$, both $\xi$ and $\eta$ can be extended as a $C^2$ 
positive function near the origin. Since the equation satisfied by $u$ is
 conformally invariant, we have
\[
F(A^{\xi})=F(A^{\eta})=1,\quad A^{\xi},~A^{\eta}\in U,\quad\mbox{in an open neighborhood of the origin}.
\]
We also know that $\xi\ge\eta$ in an open neighborhood of the origin. If (\ref{32}) does not hold, 
then $\xi (0)=\eta(0)$. By the 
arguments
 in the proof of lemma~2.1 in 
\cite{LL} which are
based on the strong maximum principle while using only
the fairly weak ellipticity hypotheses (\ref{2}) and (\ref{4}),
 we have $\xi\equiv\eta$ near the origin, 
i.e., $u(y)\equiv u_{\bar x,\bar\lambda(\bar x)}(y)$ for large $|y|$.
 Again, by the same arguments, 
$u\equiv u_{\bar x,\bar\lambda(\bar x)}$, and in particular $u$ can be extended as a $C^2$ function near the origin, 
violating our assumption that $u$ does not have such an extension. We have proved (\ref{32}).\newline
Similarly, also using arguments in the proof of lemma~2.1 in \cite{LL} 
(based on the Hopf lemma and the strong maximum principle), we have 
\begin{equation}
\frac{d}{dr}(u-u_{\bar x,\bar\lambda(\bar x)})|_{\partial B_{\bar \lambda(\bar x)}(\bar x)}>0,
\label{33}
\end{equation} 
where $\frac{d}{dr}$ denotes the outer normal differentiation with respect to $B_{\bar\lambda(\bar x)}(\bar x)$.\newline
Again, by using the strong maximum principle as in the proof of lemma~2.1 in 
\cite{LL} 
(recall that we always assume that $u$ can not be extended as a $C^1$ function near the origin), we have
\begin{equation}
(u-u_{\bar x,\bar\lambda(\bar x)})(y)>0,\quad\forall~|y-\bar x|>\bar\lambda(\bar x),~y\neq 0.
\label{34}
\end{equation}
Because of (\ref{32}), (\ref{33}),
 and the definition of $\bar\lambda(\bar x)$,   we must have, as usual, 
\begin{equation}
\liminf\limits_{y\to 0}(u-u_{\bar x,\bar\lambda (\bar x)})(y)=0.
\label{35}
\end{equation}
On the other hand, applying Lemma~\ref{lemma2} to $u$ with $w=u_{\bar x,\bar\lambda(\bar x)}$ 
(note that $\nabla u_{\bar x,\bar\lambda (\bar x)}(0)\neq\nabla v(0)$
 due to  (\ref{26}) and Lemma~\ref{lemma4}), 
we have $\liminf\limits_{x\to 0}(u-u_{\bar x,\bar\lambda(\bar x)})(x)>0$, violating (\ref{35}).
 Case 1 is settled.

{\bf In Case 2}, we have, $\forall~x\in\Bbb{R}^n\setminus\{0\}$,
\begin{equation}
u_{x,\lambda}(y)\le u(y),\quad\forall~|y-x|\ge\lambda,~y\neq 0,~0<\lambda<|x|.
\label{36}
\end{equation}
For $e\in\Bbb{R}^n$ with $\|e\|=1$ and $\mu>0$, let 
\[
\Sigma_{\mu}(e):=\{y\in\Bbb{R}^n|~y\cdot e<\mu\},\quad u^{e,\mu}(y):=u(y^{e,\mu}),
\]
where $y^{e,\mu}$ denotes the mirror symmetry point of $y$ with respect to the plane $\partial \Sigma_{\mu}(e)$.
\begin{lem}
$\forall~e\in\Bbb{R}^n$ with $\|e\|=1$ and $\forall~\mu>0$, we have
\[
u^{e,\mu}(y)\le u(y),\quad\forall~y\in\Sigma_{\mu}(e)\setminus\{0\}.
\]
\label{lemma5}
\end{lem}
{\bf Proof of Lemma~\ref{lemma5}.}\ Without loss of generality, we may assume $e=e_1$. 
For any fixed $\mu>0$, let $x=x(R)=Re_1$ for $R>\mu$, and let $\lambda=\lambda(R)=R-\mu$. By (\ref{36}),
\[
u_{x,\lambda}(y)\le u(y),\quad\forall~y\in\Sigma_{\mu}(e_1)\setminus\{0\}.
\]
Fix $y\in\Sigma_{\mu}(e_1)$, we deduce from the above that 
\[
u(y)\ge\lim\limits_{R\to\infty}u_{x,\lambda}(y)
=\lim\limits_{R\to\infty}(\frac{\lambda}{|y-x|})^{n-2}u(x+\frac{\lambda^2(y-x)}{|y-x|^2})
=u(y^{e_1,\mu}).
\]
Here we have used the fact that
 $\lim\limits_{R\to\infty}(x+\frac{\lambda^2(y-x)}{|y-x|^2})
=y^{e_1,\mu}$. 
Lemma~\ref{lemma5} is established.

\vskip 5pt
\hfill $\Box$
\vskip 5pt

It follows from Lemma~\ref{lemma5} that $w$ is radially symmetric,
 and as usual, by the Hopf Lemma 
(as in the proof of lemma~2.1 in \cite{LL}, 
using only the fairly weak ellipticity
 hypotheses (\ref{2}) and (\ref{4})), 
we have $u'(r)<0$ for $\forall~r>0$. Proposition~\ref{prop1} is established.

\vskip 5pt
\hfill $\Box$
\vskip 5pt

\begin{prop}
For $n\ge 3$, let $U\subset{\cal S}^{n\times n}$ be open and satisfy (\ref{1}) and (\ref{2}) 
and let $F\in C^1(U)$ satisfy (\ref{3}) and (\ref{4}). Assume that $u\in C^2(\Bbb{R}^n\setminus\{0\})$ 
is a positive radially symmetric function satisfying (\ref{22}), (\ref{24}) and
\begin{equation}
 u'(r)\le 0,\quad\forall~0<r<\infty.
 \label{27c}
 \end{equation}
Then either $u(r)\equiv\frac{constant}{|r|^{n-2}}$ or $u$ is of the form (\ref{7}) 
with $\bar x=0$ and some positive constants $a$ and $b$ satisfying $2b^2a^{-2}I\in U$ 
and $F(2b^2a^{-2}I)=1$.
\label{proposition2}
\end{prop}
{\bf Proof of Proposition~\ref{proposition2}.}\ If we know $\lim\limits_{r\to 0^+}(r|u'(r)|)=0$, 
then, by theorem~1.2 in \cite{LL}, $u$ is of the form (\ref{7}). By the radial symmetry of $u$, 
$\bar x=0$. Since $\infty$ is regular point of $u$, $b$ must be positive. 
Proposition~\ref{proposition2} is proved in this case. In the following, we assume that  
\begin{equation}
\limsup\limits_{r\to 0}(-ru'(r))=\limsup\limits_{r\to 0}(r|u'(r)|)>\delta>0,
\label{37}
\end{equation}
and we will show that $u(r)\equiv\frac{constant}{|r|^{n-2}}$. \newline
By (\ref{37}), we can find $r_i\to 0^+$ such that 
\begin{equation}
-r_iu'(r_i)\ge \delta,\quad\forall~i.
\label{38}
\end{equation}
Since $u$ is positive in $\Bbb{R}^n\setminus\{0\}$ and $u'(r)\le 0$ for $\forall~r>0$, we have 
$\inf\limits_{0<r<1}u(r)\ge u(1)>0$. By (\ref{24}), $\infty$ 
is a regular point of $u$. 
As usual we have, for large $\lambda>0$, that
\[
u_{\lambda}(x):=(\frac{\lambda}{|x|})^{n-2}u(\frac{\lambda^2 x}{|x|^2})\le u(x),\quad\forall~0<|x|\le \lambda.
\]
Here and below we have abused notation slightly by 
writing $u(x)=u(|x|)$.

For any fixed $i$, set 
$$\bar\lambda_i:=\{\mu>r_i~|~u_{\lambda}(x)\le u(x),~forall~r_i\le
 |x|\le\lambda,~\forall~\lambda\ge\mu\}.
$$
\begin{lem}
$\lim\limits_{i\to\infty} \bar\lambda_i=0$.
\label{lemma6}
\end{lem}
{\bf Proof of Lemma~\ref{lemma6}.}\ Suppose not, then for some positive constant $\delta_1>0$ and 
along a subsequence, we have $\bar\lambda_i>\delta>r_i$. By 
the usual arguments based on the strong maximum 
principle, the Hopf lemma and our ellipticity hypothesis, a touching must occur at $r=r_i$, i.e., 
$u_{\bar\lambda_i}(r_i)=u(r_i)$. Recall that $u_{\bar\lambda_i}(r)\le u(r)$ for $\forall~r_i\le r<\bar\lambda_i$. Thus 
\begin{equation}
u'(r_i)\ge u'_{\bar\lambda_i}(r_i).
\label{39}
\end{equation}
Since $u$ is regular at $\infty$ $\Big((\ref{24})\Big)$ and $\bar\lambda_i\ge\delta_1>0$, we have
\begin{equation}
|u'_{\bar\lambda_i}(r_i)|\le C
\label{40}
\end{equation}
for some constant $C>0$ independent of $i$. On the other hand, we have, by (\ref{38}),
\begin{equation}
\lim\limits_{i\to\infty}u'(r_i)=-\infty.
\label{41}
\end{equation}
We reach a contradiction from (\ref{39}), (\ref{40}) and (\ref{41}).
 Lemma~\ref{lemma6} is established.

\vskip 5pt
\hfill $\Box$
\vskip 5pt

\begin{lem}
$\lim\limits_{r\to 0^+} u(r)=\infty.$
\label{lemma7}
\end{lem}
{\bf Proof of Lemma~\ref{lemma7}.}\ For any fixed $\lambda>0$, we have, by 
Lemma~\ref{lemma6}, 
$\bar\lambda_i<\lambda$ for large $i$. By the definition of $\bar\lambda_i$, we have, for large $i$,
\[
u_{\lambda}(x)\le u(x),\quad\forall~r_i\le |x|\le\lambda.
\]
For any fixed $x\in \bar B_{\lambda}\setminus\{0\}$, send $i\to\infty$, 
we have $u_{\lambda}(x)\le u(x)$. 
It follows that for any  fixed $\lambda>0$, we have
\begin{eqnarray*}
\liminf\limits_{|x|\to 0}u(x)&\ge&\lim\limits_{|x|\to 0} u_{\lambda}(x)
=\lim\limits_{|x|\to 0}(\frac{\lambda}{|x|})^{n-2}u(\frac{\lambda^2 x}{|x|^2})\\
&=&\lim\limits_{|x|\to 0}\lambda^{2-n}u_{0,1}(\frac{x}{\lambda^2})=\lambda^{2-n}u_{0,1}(0).
\end{eqnarray*}
Here we have used (\ref{24}).\newline
Sending $\lambda\to 0$, we have established Lemma~\ref{lemma7}.

\vskip 5pt
\hfill $\Box$
\vskip 5pt

By Lemma~\ref{lemma7}, 
\begin{equation} 
\liminf\limits_{|x|\to\infty}(|x|^{n-2}u_{0,1}(x))=\infty.
\label{42}
\end{equation}
We also know $u_{0,1}\in C^2(\Bbb{R}^n)$ is a positive solution of 
\[
F(A^{u_{0,1}})=1,\quad A^{u_{0,1}}\in U,\quad\mbox{on}~\Bbb{R}^n.
\]
Let $w=u_{0,1}$. Starting from any point
$x\in \Bbb R^n$, the moving phere procedure can get started and can 
never stop due to (\ref{42}). 
This follows from  our usual arguments (see \cite{LZ},
\cite{LL}, \cite{LL2}).  Thus we have
\[
w_{x,\lambda}(y)\le w(y),\quad\forall
\ x\in\Bbb{R}^n, \ 0<\lambda<\infty,\
|y-x|\ge\lambda.
\] 
By a calculus lemma (see, e.g.,
lemma 11.2 in \cite{LZ}), $w\equiv$constant, i.e., 
$u(r)\equiv\frac{constant}{r^{n-2}}$. Proposition~\ref{proposition2} is established.

\vskip 5pt
\hfill $\Box$
\vskip 5pt

\noindent
{\bf Proof of Theorem~\ref{theorem1}.}\ Using the positivity and the superharmonicity of $u$ 
on $\Bbb{R}^n$, we have, by the maximum principle, 
$\liminf\limits_{|x|\to\infty}(|x|^{n-2}u(x))\ge\min\limits_{\partial B_1}u>0$. With this, we have, 
as usual, that for any
$x\in\Bbb{R}^n$, there exists some $\lambda_0(x)>0$ such that
\[
u_{x,\lambda}(y)\le u(y),\quad\forall~|y-x|\ge\lambda,~0<\lambda\le\lambda_0(x).
\]
Set, for $x\in\Bbb{R}^n$, 
$$\bar\lambda(x)
:=\{\mu>0~|~u_{x,\lambda}(y)\le u(y),~\forall~|y-x|\ge\lambda,~0<\lambda<\mu\}.
$$ 
If $\bar\lambda(x)=\infty$ for any $x\in\Bbb{R}^n$, then, as usual, $u\equiv$constant. We're done ($b=0$ in (\ref{7})). 
So, we only need to deal with the situation that $0<\bar\lambda(x)<\infty$ for some $\bar x\in\Bbb{R}^n$.\newline
The moving sphere procedure stops at $\lambda=\bar\lambda(\bar x)$,
 therefore, as usual, we have that
\begin{equation}
\liminf\limits_{|y|\to\infty}|y|^{n-2}(u-u_{\bar x,\bar\lambda(\bar x)})(y)=0,
\label{A1}
\end{equation}
\begin{equation}
(u-u_{\bar x,\bar\lambda(\bar x)})(y)\ge0,\quad\forall~|y-\bar x|\ge\bar\lambda (\bar x).
\label{A2}
\end{equation}
Let $\phi_1(x):=\bar x+\frac{\bar\lambda (\bar x)^2(x-\bar x)}{|x-\bar x|^2}$,
 we know 
$u_{\phi_1}=u_{\bar x,\bar\lambda(\bar x)}$,
where $u_{\phi_1}
:=|J_{\phi_1}|^{ \frac {n-2}{2n} }
(u\circ \phi_1)$, $J_{\phi_1}$ denotes the Jacobian of $\phi_1$.
 Pick any $\tilde x\neq\bar x$ and let
\[
\phi_2(x):=\tilde x+\frac{x-\tilde x}{|x-\tilde x|^2},\quad\tilde u:=u_{\phi_2},
\quad\tilde v:=(u_{\phi_1})_{\phi_2}=u_{\phi_1 \circ \phi_2}.
\]
Then $\tilde u\in C^2(\Bbb{R}^n\setminus\{\tilde x\})$, $\infty$ is a
 regular point of 
$\tilde u$ (i.e., $\frac{1}{|x|^{n-2}}\tilde u(\frac{x}{|x|^2})$ can be extended to a positive 
$C^2$ function near the origin), $\Delta\tilde u\le 0$ in $\Bbb{R}^n\setminus\{\tilde x\}$, 
$\tilde v\in C^2(\Bbb{R}^n\setminus\{\phi_2^{-1}(\bar x)\})$, $\infty$ is a regular point of 
$\tilde v$ (since $\bar x\neq\tilde x$), $\Delta\tilde v\le 0$ in $\Bbb{R}^n\setminus\{\phi_2^{-1}(\bar x)\}$, 
$\tilde u\ge\tilde v$ in an open
neighborhood of $\tilde x$ (because of (\ref{A2})), and 
$\liminf\limits_{x\to\tilde x}(\tilde u-\tilde v)(x)=0$ 
(because of (\ref{A1})). By (\ref{1}) and 
the conformal invariance of the equation satisfied by $u$, we have 
\[
F(A^{\tilde u})=1,\quad A^{\tilde u}\in U,\quad\mbox{in}~\Bbb{R}^n\setminus\{\tilde x\}.
\]
Since $\tilde x\neq\bar x$, we have $\phi_2^{-1}(\bar x)\neq\tilde x$, therefore $\tilde v$ 
is a positive $C^2$ function near $\tilde x$. If $\nabla\tilde v(\tilde x)=0$, then, by 
applying 
Proposition~\ref{prop1} to $\hat u(x):=\tilde u(\tilde x+x)$, $\hat u$ is radially symmetric and
\[
\hat u'(r)<0,\quad\forall~0<r<\infty.
\]
Next, by applying Proposition~\ref{proposition2} to $\hat u$, we have either
\begin{equation}
\hat u(x)\equiv\frac{constant}{|x|^{n-2}},
\label{A3}
\end{equation}
or, for some positive constants $a$ and $b$,
\begin{equation}
\hat u(r)\equiv (\frac{a}{1+b^2r^2})^{\frac{n-2}{2}}.
\label{A4}
\end{equation}
If (\ref{A3}) occurs, then $u\equiv$constant, i.e., $u$ is of the form (\ref{7}) with 
$b=0$ and some $a>0$. If (\ref{A4}) occurs, then
\[
u(y)=\frac{1}{|y-\tilde x|^{n-2}}\hat u(\frac{1}{|y-\tilde x|})=(\frac{a}{b^2+|y-\tilde x|^2})^{\frac{n-2}{2}},
\]
and therefore $u$ is of the form (\ref{7}). Thus we have proved 
Theorem~\ref{theorem1} provided that 
$\nabla\tilde v(\tilde x)=0$. If $\nabla\tilde v(\tilde x)\neq 0$, 
we will make a suitable M\"obius 
transformation to reduce it
 to the situation with $\nabla\tilde v(\tilde x)=0$. For this, we need the following 
fact (used in the proof of theoerm~1.1 in \cite{LL}).
\begin{lem}
Let $s>0$, $y,p\in\Bbb{R}^n\setminus\{0\}$ with $n\ge 3$ and $y=\frac{(2-n)s}{|p|^2}p$. Assume that 
$\xi$ is a $C^1$ function near $y$ satisfying $\xi (y)=s$ and $\nabla\xi(y)=p$. Then
\[
(\nabla\xi_{\psi})(\psi^{-1}(y))=0,
\]
where $\psi(x):=\frac{\lambda^2x}{|x|^2}$ for any fixed $\lambda>0$.
\label{lemma8}
\end{lem}
Lemma~\ref{lemma8} follows from a direct computation.

Back to the proof of Theorem~\ref{theorem1}, when $\nabla\tilde v(\tilde x)\neq 0$, let 
$s=\tilde v(\tilde x)>0$, $p=\nabla \tilde v(\tilde x)\neq 0$, and $y=\frac{(2-n)s}{|p|^2}p$. 
Define $\xi(x):=\tilde v(x-y+\tilde x)$, $\psi(x):=\frac{|y|^2x}{|x|^2}$.
 By Lemma~\ref{lemma8}, 
$(\nabla\xi_{\psi})(\psi^{-1}(y))=0$. Now let 
\[
\eta(x)=\tilde u(x-y+\tilde x),\quad\hat u=\eta_{\psi},\quad\hat v=\xi_{\psi}.
\]
Then $\hat u\in C^2(\Bbb{R}^n\setminus\{\psi^{-1}(y)\})$, $\infty$ is a regular point of $\hat u$, 
$\Delta \hat u\le 0$ in $\Bbb{R}^n\setminus\{\psi^{-1}(y)\}$, $\hat v$ is a positive 
$C^2$ superharmonic function in an open neighborhood of $\psi^{-1}(y)$, 
$\hat u\ge\hat v$ in an open neighborhood of $\psi^{-1}(y)$, 
\[
\liminf\limits_{x\to\psi^{-1}(y)}(\hat u-\hat v)(x)=0,
\]
and
\[
F(A^{\hat u})=1,\quad A^{\hat u}\in U,\quad\mbox{in}~\Bbb{R}^n\setminus\{\psi^{-1}(y)\}.
\]
Now we also know that $\nabla\hat v(\psi^{-1}(y))=0$. So we have, by applying 
Proposition~\ref{prop1} to $u^{\ast}(x):=\hat u(x+\psi^{-1}(y))$, that $u^{\ast}$ is radially symmetric and 
\[
(u^{\ast})'(r)<0,\quad\forall~0<r<\infty.
\]
Applying Proposition~\ref{proposition2} to $u^{\ast}$, we have either
\begin{equation}
u^{\ast}(x)\equiv\frac{constant}{|x|^{n-2}},
\label{A5}
\end{equation}
or, for some positive constants $a$ and $b$,
\begin{equation}
u^{\ast}(r)\equiv(\frac{a}{1+b^2r^2})^{\frac{n-2}{2}}.
\label{A6}
\end{equation}
If (\ref{A5}) occurs, we have $u\equiv$constant. If (\ref{A6}) occurs, 
$u$ is of the form (\ref{7}) and $u$ is not a constant. Theorem~\ref{theorem1} is established.

\vskip 5pt
\hfill $\Box$

\end{document}